\documentclass[12pt]{article}
\usepackage{amsfonts,amscd,a4}
\usepackage{color}
\usepackage{graphicx}
\usepackage{theorem}
\newtheorem{theo}{Theorem}
\newtheorem{prop}[theo]{Proposition}

\newtheorem{coro}[theo]{Corollary}
\newtheorem{defi}[theo]{Definition}
{\theorembodyfont{\rm}
\newtheorem{example}[theo]{Example}
}
\newcommand{\bA}{{\bf A}}
\newcommand{\bB}{{\bf B}}

\newcommand{\bI}{{\bf I}}
\newcommand{\bJ}{{\bf J}}
\newcommand{\bK}{{\bf K}}

\newcommand{\cA}{{\mathcal A}}
\newcommand{\cB}{{\mathcal B}}
\newcommand{\cC}{{\mathcal C}}

\newcommand{\cF}{{\mathcal F}}
\newcommand{\cG}{{\mathcal G}}

\newcommand{\cI}{{\mathcal I}}
\newcommand{\cJ}{{\mathcal J}}
\newcommand{\cK}{{\mathcal K}}
\newcommand{\cL}{{\mathcal L}}

\newcommand{\cP}{{\mathcal P}}

\newcommand{\cS}{{\mathcal S}}

\newcommand{\eT}{{\sf T}}

\newcommand{\sC}{{\mathbb C}}

\newcommand{\sM}{{\mathbb M}}
\newcommand{\sN}{{\mathbb N}}

\newcommand{\sR}{{\mathbb R}}

\newcommand{\sT}{{\mathbb T}}

\newcommand{\sZ}{{\mathbb Z}}
\newcommand{\qed}{\rule{1ex}{1ex}}

\newcommand{\diag}{\mbox{\rm diag} \,}

\newcommand{\ess}{\mbox{\rm ess} \,}

\newcommand{\im}{{\rm im} \,}

\newcommand{\rank}{\mbox{\rm rank} \,}

\begin{document}
\title{Arveson dichotomy and essential fractality}
\author{Steffen Roch}
\date{Dedicated to Vladimir S. Rabinovich on the occasion of his 70th birthday}
\maketitle
\begin{abstract}
The notions of fractal and essentially fractal algebras of approximation sequences and of the Arveson dichotomy have proved extremely useful for several spectral approximation problems. The purpose of this short note is threefold: to present a short new proof of the fractal restriction theorem, to relate essential fractality with Arveson dichotomy, and to derive a restriction theorem for essential fractality.
\end{abstract}
{\bf 2010 AMS-MSC:} 65J10, 46L99, 47N40
\section{Preliminaries} \label{s1}
Let $H$ be an infinite dimensional separable Hilbert space. We denote by $L(H)$ the $C^*$-algebra of the bounded linear operators and by $K(H)$ the ideal of the compact operators on $H$.

A sequence  $\cP = (P_n)_{n \ge 1}$ of orthogonal projections of finite rank which converge strongly to the identity operator on $H$ is called a {\em filtration} on $H$. Given a filtration $\cP$, let $\cF^\cP$ stand for the set of all sequences $\bA = (A_n)$ of operators $A_n : \im P_n \to \im P_n$ such that the sequence $(A_n P_n)$ converges strongly to an operator $W^\cP(\bA) \in L(H)$. Since every sequence in $\cF^\cP$ is bounded by the Banach-Steinhaus theorem, one can introduce pointwise defined operations
\begin{equation} \label{e12.1}
(A_n) + (B_n) := (A_n + B_n), \quad (A_n)(B_n) := (A_n B_n), \quad (A_n)^* := (A_n^*)
\end{equation}
and the supremum norm $\|(A_n)\|_\cF := \sup_n \|A_n\|$, which make $\cF^\cP$ to a unital $C^*$-algebra and $W^\cP : \cF^\cP \to L(H)$ to a unital $^*$-homomorphism. This homomorphism is also known as the consistency map associated with the filtration $\cP$.

Set $\delta(n) := \rank P_n := \dim \im P_n < \infty$ for every $n$ and choose an orthonormal basis in each of the spaces $\im P_n$. Every operator $A_n \in L(\im P_n)$ can be identified with its matrix representation with respect to the chosen basis and, thus, with an element of the $C^*$-algebra $\sC^{\delta(n) \times \delta(n)}$ of all $\delta(n) \times \delta(n)$ matrices with complex entries. The choice of a basis in each space $\im P_n$ makes $\cF^\cP$ to a special instance of an algebra of matrix sequences in the following sense. Given a sequence $\delta$ of positive integers, we let $\cF^\delta$ stand for the set of all bounded sequences $(A_n)$ of matrices $A_n \in \sC^{\delta(n) \times \delta(n)}$. Introducing again pointwise operations and the supremum norm, we make $\cF^\delta$ to a $C^*$-algebra with identity element $(I_{\delta(n)})$, the {\em algebra of matrix sequences with dimension function} $\delta$. The set of all sequences in $\cF^\delta$ which tend to zero in the norm forms a closed ideal of $\cF^\delta$ which we denote by $\cG^\delta$. For example, the algebra of matrix sequences with constant dimension function $\delta = 1$ is $l^\infty(\sN)$, but in what follows we will be mainly interested in strictly increasing dimension functions, as they occur in the context of filtrations.

When passing from $\cF^\cP$ to $\cF^\delta$ with $\delta(n) := \rank P_n$, one loses the embedding of the matrix algebras $L(\im P_n) \cong \sC^{\delta(n) \times \delta(n)}$ into a common Hilbert space. It makes thus no sense to speak about strong convergence of a sequence in $\cF^\delta$. But it will turn out that algebras of matrix sequences provide a suitable frame to formulate and study stability problems as well as a lot of other problems which do not depend upon an embedding into a Hilbert space. Moreover, some of the notions and assertions discussed in this paper remain meaningful in the much more general context, when $\cF^\cC$ is the direct product of a sequence $\cC = (\cC_n)_{n \ge 1}$ of unital $C^*$-algebras. The associated ideal of zero sequences in $\cF^\cC$, which can be identified with the direct sum of the family $\cC$ in a natural way, will then be denoted by $\cG^\cC$.

The following will serve as a running example in this paper. We consider the algebra of the finite sections discretization for Toeplitz operators with continuous generation function. For a continuous function $a$ on the complex unit circle $\sT$, the associated {\em Toeplitz operator} is the operator $T(a)$ on $l^2(\sZ^+)$ which is given by the infinite matrix $(a_{i-j})_{i,j = 0}^\infty$, with $a_k$ denoting the $k$th Fourier coefficient of $a$. Note that $T(a)$ is a bounded operator and $\|T(a)\| = \|a\|_\infty$. For $n \in \sN$, put
\[
P_n : l^2(\sZ^+) \to l^2(\sZ^+), \quad (x_n)_{n \ge 0} \mapsto (x_0, \, x_1, \, \ldots, \, x_{n-1}, \, 0, \, 0, \, \ldots).
\]
Then $\cP = (P_n)$ is a filtration on $l^2(\sZ^+)$. We let $\cS(\eT(C))$ stand
for the smallest closed subalgebra of $\cF^\cP$ which contains all sequences
$(P_n T(a)|_{\im P_n})$ of finite sections of Toeplitz operators
$T(a)$ with $a \in C(\sT)$. Let $R_n : \im P_n \to \im P_n$ be the reflection
operator
\[
(x_0, \, x_1, \, \ldots, \, x_{n-1}, \, 0, \, 0, \, \ldots) \mapsto (x_{n-1}, \, \ldots, \, x_1, \, x_0, \, 0, \, 0, \, \ldots).
\]
It is not hard to see that for each sequence $\bA = (A_n) \in \cS(\eT(C))$, the strong limit $\widetilde{W} (\bA) :=
\mbox{s-lim} R_n A_n R_n P_n$ exists and that $\widetilde{W}$ is a unital and fractal $^*$-homomorphism from $\cS(\eT(C))$ to $L(l^2\sZ^+)$. The following is a by now classical result by B\"ottcher and Silbermann \cite{BSi8}, see also Chapter 2 in \cite{BSi2} and Sections 1.3, 1.4 and 1.6 in \cite{HRS2}. 
\begin{theo} \label{tneu}
$(a)$ The algebra $\cS(\eT(C))$ consists of all sequences $(P_n T(a) P_n + P_n K P_n + R_nLR_n + G_n)$ where $a \in C(\sT)$, $K, \, L \in K(l^2(\sZ^+))$, and $(G_n) \in \cG^\cP$. \\[1mm]
$(b)$ For every sequence $\bA \in \cS(\eT(C))$, the coset $\bA + \cG^\cP$ is invertible in the quotient algebra $\cS(\eT(C))/\cG^\cP$ if and only if the operators $W^\cP(\bA)$ and $\widetilde{W} (\bA)$ are invertible.
\end{theo}
Due to its transparent structure, the algebra $\cS(\eT(C))$ served as a basic example for the development of algebraic methods in asymptotic numerical analysis. These methods have found fruitful applications in the stability analysis of different approximation methods for numerous classes of operators; see the monographs \cite{HRS2,PrS1,RSS2} for an overview. In particular, I would like to emphasize the finite sections method for band-dominated operators, a topic which was mainly influenced and shaped by Vladimir S. Rabinovich and the limit operator techniques developed by him, see \cite{RaR5,RRS1,RRS2,RRSB,RRS6} and \cite{Roc10} for an overview. In fact, the algebra of the finite sections method for band-dominated operators is the first real-life example of an essentially fractal, but not fractal, algebra (these notions will be introduced below).
\section{Fractality} \label{ssf1.1a}
As it was observed in \cite{Roc0,RoS5}, several natural approximation procedures lead to $C^*$-subalgebras $\cA$ of the algebra $\cF$ which are distinguished by the property of self-similarity: Given a subsequence of a sequence in $\cA$, one can uniquely reconstruct the full sequence up to a sequence which tends to zero in the norm. These algebras were called fractal in \cite{RoS5}. The goals of this section is to recall the basic definitions and some consequences of fractality, and to give a short proof of the known fact that every separable subalgebra of $\cF$ possesses a fractal restriction.

In this section, we let $\cF := \cF^\cC$ be the product of a family $\cC = (\cC_n)_{n \in \sN}$ of unital $C^*$-algebras and $\cG := \cG^\cC$ the associated ideal of zero sequences.
\subsection{Definition and first consequences}
For each strictly increasing sequence $\eta :\sN \to \sN$, let $\cF_\eta$ stand for the product of the family $(\cC_{\eta(n)})_{n \in \sN}$ of $C^*$-algebras, and write $\cG_\eta$ for the associated ideal of zero sequences. The elements of $\cF_\eta$ can be viewed of as subsequences of sequences in $\cF$. The canonical restriction mapping $R_\eta: \cF \to \cF_\eta, \; (A_n) \mapsto (A_{\eta(n)})$ is a $^*$-homomorphism from $\cF$ onto $\cF_\eta$ and maps $\cG$ onto $\cG_\eta$. More generally, for each $C^*$-subalgebra $\cA$ of $\cF$, we let $\cA_\eta$ denote the image of $\cA$ under $R_\eta$. Clearly, $\cA_\eta$ is a $C^*$-subalgebra of $\cF_\eta$. We call algebras obtained in this way {\em restrictions} of $\cA$.
\begin{defi} \label{df1.2}
$(a)$ Let $\cA$ be a $C^*$-subalgebra of $\cF$. A $^*$-homomorphism $W$ from $\cA$ into a $C^*$-algebra $\cB$ is called {\em fractal} if it factors through $R_\eta|_\cA$ for {\em every} strictly increasing sequence $\eta : \sN \to \sN$, i.e., if for each such $\eta$, there is a mapping $W_\eta : \cA_\eta \to \cB$ such that $W = W_\eta R_\eta|_\cA$. \\[1mm]
$(b)$ A $C^*$-subalgebra $\cA$ of $\cF$ is {\em fractal} if the canonical homomorphism
\[
\cA \to \cA/(\cA \cap \cG), \quad \bA \mapsto \bA + (\cA \cap \cG)
\]
is fractal. \\[1mm]
$(c)$ A sequence $\bA \in \cF$ is {\em fractal} if the smallest $C^*$-subalgebra of $\cF$ which contains the sequence $\bA$ and the identity sequence is fractal.
\end{defi}
For example, if $\cP$ is a filtration, then the associated consistency map $W^\cP$ is fractal (since the strong limit of a sequence $(A_n) \in \cF^\cP$ can be determined from each subsequence of $(A_n)$). For the same reason, the homomorphism $\widetilde{W}$ appearing in Theorem \ref{tneu} is fractal.

The fractal subalgebras of $\cF$ are distinguished by their property that every sequence in the algebra can be rediscovered from each of its (infinite) subsequences up to a sequence tending to zero. Note that, by Definition \ref{df1.2}, a fractal sequence always lies in a {\em unital} fractal algebra, whereas a fractal algebra needs not to be unital.

Assertion $(a)$ of the following theorem provides an equivalent characterization of the fractality of an algebra. Proofs of Theorems \ref{tf1.8} and \ref{tf1.12} are given in \cite{RoS5} and in Section 1.6 of \cite{HRS2}.
\begin{theo} \label{tf1.8}
$(a)$ A $C^*$-subalgebra $\cA$ of $\cF$ is fractal if and only if the implication
\begin{equation} \label{ef1.9}
R_\eta (\bA) \in \cG_\eta \; \Rightarrow \; \bA \in \cG
\end{equation}
holds for every sequence $\bA \in \cA$ and every strictly increasing sequence $\eta$. \\[1mm]
$(b)$ If $\cA$ is a fractal $C^*$-subalgebra of $\cF$, then $\cA_\eta \cap \cG_\eta = (\cA \cap \cG)_\eta$ for each strictly increasing sequence $\eta$. \\[1mm]
$(c)$  A unital $C^*$-subalgebra of $\cF$ is fractal if and only if each of its elements is fractal.
\end{theo}
The following criterion will prove to be useful in order to verify the fractality of many specific algebras of approximation methods.
\begin{theo} \label{tf1.12}
A unital $C^*$-subalgebra $\cA$ of $\cF$ is fractal if and only if there is a family $\{W_t\}_{t \in T}$ of unital and fractal $^*$-homomorphisms $W_t$ from $\cA$ into unital $C^*$-algebras $\cB_t$ such that the following equivalence holds for every sequence $\bA \in \cA$: The coset $\bA + \cA \cap \cG$ is invertible in $\cA/(\cA \cap \cG)$ if and only if $W_t(\bA)$ is invertible in $\cB_t$ for every $t \in T$.
\end{theo}
For example, since $W^\cP$ and $\widetilde{W}$ are fractal homomorphisms, we conclude from Theorem \ref{tneu} $(b)$ and from the previous theorem that the algebra $\cS(\eT(C))$ is fractal. \hfill \qed \\[3mm]
The property of fractality has striking consequences for asymptotic spectral properties of a sequence $\bA = (A_n)$, see \cite{Roc0,RoS5} and Chapter 3 in \cite{HRS2}. Here we only mention a few of them which are relevant for what follows. For every element $a$ of a unital $C^*$-algebra $\cA$, we let $\sigma_2(a)$ denote the set of all non-negative square roots of points in the spectrum of $a^*a$. In case $\cA = \sC^{n \times n}$, the numbers in $\sigma_2(a)$ are known as the singular values of $a$.
\begin{prop} \label{pf2.1}
Let $\cA$ be a fractal $C^*$-subalgebra of $\cF$ and $\bA = (A_n)$ a sequence in $\cA$. Then \\[1mm]
$(a)$ the sequence $\bA$ is stable if and only if it possesses a stable subsequence; \\[1mm]
$(b)$ the limit $\lim_{n \to \infty} \|A_n\|$ exists and is equal to $\|\bA + \cG\|$; \\[1mm]
$(c)$ the limit $\lim_{n \to \infty} \sigma_2(A_n)$ exists with respect to the Hausdorff distance on $\sR$ and is equal to $\sigma_2(\bA + \cG)$.
\end{prop}
\subsection{The fractal restriction theorem} \label{ssf3.0}
The preceding proposition and related results from \cite{HRS2} indicate that it is a question of vital importance in numerical analysis to single out fractal subsequences of a given sequence in $\cF$. The following theorem states that such subsequences always exist.
\begin{theo} \label{tf3.20}
Let $\cA$ be a separable $C^*$-subalgebra of $\cF$. Then there exists a strictly increasing sequence $\eta : \sN \to \sN$ such that the restricted algebra $\cA_\eta = R_\eta \cA$ is a fractal subalgebra of $\cF_\eta$.
\end{theo}
Since finitely generated $C^*$-algebras are separable, this result immediately implies:
\begin{coro} \label{cf3.21}
Every sequence in $\cF$ possesses a fractal subsequence.
\end{coro}
Theorem \ref{tf3.20} was first proved in \cite{Roc0}. We shall give a much shorter proof here, which is based on the following converse of assertion $(b)$ of Proposition \ref{pf2.1} (whereas the original proof used the converse of assertion $(c)$ of this proposition).
\begin{prop} \label{pf3.0a}
Let $\cA$ be a $C^*$-subalgebra of $\cF$ and $\cL$ a dense subset of $\cA$. If the sequence of the norms $\|A_n\|$ converges for each sequence $(A_n) \in \cL$, then the algebra $\cA$ is fractal.
\end{prop}
{\bf Proof.} First we show that if the sequence of the norms converges for each sequence in $\cL$, then it converges for each sequence in $\cA$. Let $(A_n) \in \cA$ and $\varepsilon > 0$. Choose $(L_n) \in \cL$ such that $\|(A_n - L_n)\| = \sup \|A_n - L_n\| < \varepsilon/3$, and let $n_0 \in \sN$ be such that $\left| \|L_n\| - \|L_m\| \right| < \varepsilon$ for all $m, \, n \ge n_0$. Then, for $m, \, n \ge n_0$,
\begin{eqnarray*}
\left| \|A_n\| - \|A_m\| \right| & \le &  \left| \|A_n\| - \|L_n\| \right| + \left| \|L_n\| - \|L_m\| \right| + \left| \|L_m\| - \|A_m\| \right| \\
& \le & \|A_n - L_n\| + \left| \|L_n\| - \|L_m\| \right| + \|L_m - A_m\| \; \le \; \varepsilon.
\end{eqnarray*}
Thus, $(\|A_n\|)$ is a Cauchy sequence, hence convergent. But the convergence of the norms for each sequence in $\cA$ implies the fractality of $\cA$ by Theorem \ref{tf1.8}. Indeed, if a subsequence of a sequence $(A_n) \in \cA$ tends to zero, then $0 = \liminf \|A_n\| = \lim \|A_n\|$, whence $(A_n) \in \cG$. \hfill \qed \\[3mm]
{\bf Proof of Theorem \ref{tf3.20}.} Let $\{\bA^m\}_{m \in \sN}$ of $\cA_{sa}$ be a dense countable subset of $\cA$ which consists of sequences $\bA^m = (A_n^m)_{n \in \sN}$. Let $\eta_1 : \sN \to \sN$ be a strictly increasing sequence such that the sequence of the norms $\|A_{\eta_1(n)}^1\|$ converges. Next let $\eta_2$ be a strictly increasing subsequence of $\eta_1$ such that the sequence $(\|A_{\eta_2(n)}^2\|)_{n \in \sN}$ converges. We proceed in this way and find, for each $k \ge 2$, a strictly increasing subsequence $\eta_k$ of $\eta_{k-1}$ such that the sequence $(\|A_{\eta_k(n)}^k\|)_{n \in \sN}$ converges. Define the sequence $\eta$ by $\eta(n) := \eta_n(n)$. Then $\eta$ is  strictly increasing, and the sequence $(\|A_{\eta(n)}^k\|)_{n \in \sN}$ converges for every $k \in \sN$.

Since the sequences $R_\eta (\bA^m)$ with $k \in \sN$ form a dense subset of the restricted algebra $\cA_\eta$, and since each sequence $R_\eta (\bA^m) = (A_{\eta(n)}^k)_{n \in \sN}$ has the property that the sequence of the norms $\|A_{\eta(n)}^k \|$ converges, the assertion follows from Proposition \ref{pf3.0a}. \hfill \qed
\section{Essential fractality} \label{s43}
Recall that a $C^*$-subalgebra $\cA$ of $\cF$ is fractal if each sequence $(A_n) \in \cA$ can be rediscovered from each of its (infinite) subsequences modulo a sequence in the ideal $\cG$. There are plenty of subalgebras of $\cF$ which arise from concrete discretization methods and which are fractal (the finite sections algebra $\cS(\eT(C))$ for Toeplitz operators is one example). On the other hand, the algebra of the finite sections method for band-dominated operators is an example of an algebra which fails to be fractal. But the latter algebra enjoys a weaker form of fractality which we called essential fractality in \cite{Roc10}. Basically, a $C^*$-subalgebra $\cA$ of $\cF$ is essentially fractal if each sequence $(A_n) \in \cA$ can be rediscovered from each of its (infinite) subsequences modulo a sequence in the ideal $\cK$ of the compact sequences. The role of this ideal in numerical analysis can be compared with the role of the ideal of the compact operators in operator theory.

In this section, we first recall the definition of a compact sequence and state some useful characterizations of compactness and the definitions of $\cJ$-fractality and essential fractality from \cite{Roc10}. The main goal of this section is to derive an analogue of the fractal restriction theorem for essential fractality.

Unless otherwise stated, we let $\cF = \cF^\delta$ be an algebra of matrix sequences with dimension function $\delta$ and $\cG := \cG^\delta$ the associated ideal of zero sequences in this section.
\subsection{Compact sequences} \label{s41}
Slightly abusing the notation, we call a sequence $(K_n) \in \cF$ a {\em sequence of rank one matrices} if the rank of every matrix $K_n$ is less than or equal to one. The product of a sequence of rank one matrices with a sequence in $\cF$ is a sequence of rank one matrices again. Hence, the set of all finite sums of sequences of rank one matrices forms an (in general, non-closed) ideal of $\cF$. We let $\cK$ denote the closure of this ideal and refer to the elements of $\cK$ as {\em compact sequences}. Thus, $\cK$ is the smallest closed ideal of $\cF$ which contains all sequences of rank one matrices, and a sequence $(A_n) \in \cF$ is compact if, and only if, for every $\varepsilon > 0$, there is a sequence $(K_n) \in \cF$ such that
\begin{equation} \label{e41.1}
\sup_n \|A_n - K_n\| < \varepsilon \quad \mbox{and} \quad \sup_n \, \rank K_n < \infty.
\end{equation}
Note that $\cK$ contains the ideal $\cG$, and that the restriction of a compact sequence is compact. More precisely, if $\bK$ is a compact sequence in the algebra $\cF^\delta$ of matrix sequences with dimension function $\delta$ and if $\eta$ is a strictly increasing sequence, then the restriction $R_\eta \bK$ is a compact sequence in the algebra $R_\eta \cF^\delta \cong \cF^{\delta \circ \eta}$ of matrix sequences with dimension function $\delta \circ \eta$.

An appropriate notion of the rank of a sequence in $\cF$ can be introduced as follows. A sequence $\bA \in \cF$ has {\em finite essential rank} if it is the sum of a sequence in $\cG$ and a sequence $(K_n)$ with $\sup_n \rank K_n < \infty$. If $\bA$ is of finite essential rank, then there is a smallest integer $r \ge 0$ such that $\bA$ can be written as $(G_n) + (K_n)$ with $(G_n) \in \cG$ and $\sup_n \rank K_n \le r$. We call this integer the {\em essential rank}of $\bA$ and write $\ess \rank \bA = r$. Thus, the sequences of essential rank 0 are just the sequences in $\cG$. If $\bA$ is not of finite essential rank, we set $\ess \rank \bA = \infty$. Clearly, the sequences of finite essential rank form an ideal of $\cF$ which is dense in $\cK$, and for arbitrary sequences $\bA, \bB \in \cF$ one has
\[
\ess \rank (\bA + \bB) \le \ess \rank \bA + \ess \rank \bB,
\]
\[
\ess \rank (\bA \bB) \le \min \, \{ \ess \rank \bA, \, \ess \rank \bB \}.
\]
Given a filtration $\cP = (P_n)$ on a Hilbert space $H$, we identify the algebra $\cF^\cP$ with the algebra $\cF$ of matrix sequences with dimension function $\delta(n) := \rank P_n$. Note that this identification requires the choice of an orthogonal basis in each space $\im P_n$. We define the {\em ideal $\cK^\cP$ of the compact sequences in} \index{$\cK^\cP$} $\cF^\cP$ in the same way as before. It is clear that then the ideal $\cK^\cP$ can be identified with $\cK$, independently of the choice of the bases.

For example, using the explicit description of the finite sections algebra of Toeplitz operators in Theorem \ref{tneu} $(a)$, it is not hard to show that the intersection $\cS(\eT(C)) \cap \cK$ consists of all sequences
\begin{equation} \label{eneu1}
(P_n K P_n + R_n LR_n + G_n) \quad \mbox{with} \; K, \, L \; \mbox{compact and} \; (G_n) \in \cG
\end{equation}
and that the essential rank of the sequence (\ref{eneu1}) is equal to $\rank K + \rank L$. \hfill \qed \\[3mm]
There are several equivalent characterizations of compact sequences, see \cite{Roc10}. In what follows we shall need a characterization of a compact sequence $(K_n)$ in terms of the asymptotic behavior of the singular values of the entries $K_n$. To state this criterion, we denote the decreasingly ordered singular values of an $n \times n$ matrix $A$ by
\begin{equation} \label{e41.5}
\|A\| = \Sigma_1(A) \ge \Sigma_2(A) \ge \ldots \ge  \Sigma_n (A) \ge 0
\end{equation}
and recall from Linear Algebra that $A^*A$ and $AA^*$ are unitarily equivalent, whence $\Sigma_k(A) = \Sigma_k (A^*)$, and that every matrix $A$ has a singular value decomposition (SVD)
\begin{equation} \label{e41.6a}
A = E^* \, \diag (\Sigma_1 (A), \, \ldots, \, \Sigma_n (A)) F
\end{equation}
with unitary matrices $E$ and $F$.

The announced characterization of compact sequences in terms of singular values reads as follows. See Sections 4.2 and 5.1 in \cite{Roc10} for the proof of this and the following theorem.
\begin{theo} \label{t41.12}
The following conditions are equivalent for a sequence $(K_n) \in \cF$: \\[1mm]
$(a)$ $\lim_{k \to \infty} \sup_{n \ge k} \Sigma_k (K_n) = 0$; \\[1mm]
$(b)$ $\lim_{k \to \infty} \limsup_{n \to \infty} \Sigma_k (K_n) = 0$; \\[1mm]
$(c)$ the sequence $(K_n)$ is compact.
\end{theo}
A sequence in $\cF$ is called a {\em Fredholm sequence} if it is invertible modulo $\cK$. As the compact sequences, Fredholm sequences can be characterized in terms of singular values. Let $\sigma_1(A) \le \ldots \le \sigma_n(A)$ denote the increasingly ordered singular values of an $n \times n$-matrix $A$.
\begin{theo} \label{t42.1}
The following conditions are equivalent for a sequence $(A_n) \in \cF$: \\[1mm]
$(a)$ $(A_n)$ is a Fredholm sequence. \\[1mm]
$(b)$ There are sequences $(B_n) \in \cF$ and $(J_n) \in \cK$ with $\sup_n \rank J_n < \infty$ such that
$B_n A_n = I_n + J_n$ for all $n \in \sN$. \\[1mm]
$(c)$ There is a $k \in \sZ^+$ such that $\liminf_{n \to \infty} \sigma_{k+1} (A_n) > 0$.
\end{theo}
\subsection{$\cJ$-fractal algebras} \label{ss43.1}
Our next goal is to introduce fractality of an algebra $\cA$ with respect to an arbitrary ideal $\cJ$ in place of $\cG$. The results presented in this subsection hold in the general case, when $\cF$ is the product of a family $(\cC_n)_{n \in \sN}$ of unital $C^*$-algebras. We start with a criterion for the fractality of the canonical quotient map $\cA \to \cA/\cJ$.
\begin{theo} \label{t43.1}
Let $\cA$ be a $C^*$-subalgebra of $\cF$ and $\cJ$ a closed ideal of $\cA$. The canonical homomorphism $\pi^\cJ : \cA \to \cA/\cJ$ is fractal if and only if the following implication holds for every sequence $\bA \in \cA$ and every strictly increasing sequence $\eta : \sN \to \sN$
\begin{equation} \label{e43.2}
R_\eta (\bA) \in \cJ_\eta \quad \Longrightarrow \quad \bA \in \cJ.
\end{equation}
\end{theo}
{\bf Proof.} Let $\pi^\cJ$ be fractal, i.e., for each $\eta$, there is a mapping $\pi^\cJ_\eta$ such that $\pi^\cJ = \pi^\cJ_\eta R_\eta|_\cA$. Let $R_\eta (\bA) \in \cJ_\eta$ for a sequence $\bA \in \cA$. We choose a sequence $\bJ \in \cJ$ such that $R_\eta (\bA) = R_\eta (\bJ)$. Applying the homomorphism $\pi^\cJ_\eta$ to both sides of this equality we obtain $\pi^\cJ (\bA) = \pi^\cJ (\bJ) = 0$, whence $\bA \in \cJ$.

For the reverse implication, let $\bA$ and $\bB$ be sequences in $\cA$ with $R_\eta (\bA) = R_\eta (\bB)$. Then $R_\eta (\bA - \bB) = 0 \in \cJ_\eta$, and (\ref{e43.2}) implies that $\bA - \bB \in \cJ$. Thus, the mapping
\[
\pi^\cJ_\eta : \cA_\eta \to \cA/\cJ, \quad R_\eta (\bA) \mapsto \bA + \cJ
\]
is correctly defined, and it satisfies $\pi^\cJ_\eta R_\eta|_\cA = \pi^\cJ$.  \hfill \qed \\[3mm]
Let now $\cJ$ be a closed ideal of $\cF$. Then $\cA \cap \cJ$ is a closed ideal of $\cA$, and the preceding theorem states that the canonical mapping $\pi^{\cA \cap \cJ} : \cA \to \cA/(\cA \cap \cJ)$ is fractal if and only if the implication
\begin{equation} \label{e43.2a}
R_\eta (\bA) \in (\cA \cap \cJ)_\eta \quad \Longrightarrow \quad \bA \in \cJ
\end{equation}
holds for every sequence $\bA \in \cA$ and every strictly increasing sequence $\eta$. It would be much easier to check this implication if one would have
\begin{equation} \label{e43.2b}
(\cA \cap \cJ)_\eta = \cA_\eta \cap \cJ_\eta
\end{equation}
foe every $\eta$, in which case the implication (\ref{e43.2a}) reduces to $R_\eta (\bA) \in \cJ_\eta \Rightarrow \bA \in \cJ$. Recall from Theorem \ref{tf1.8} $(b)$ that (\ref{e43.2b}) indeed holds if $\cJ = \cG$ and if the canonical homomorphism $\pi^{\cA \cap \cG} : \cA \to \cA/(\cA \cap \cG)$ is fractal. One cannot expect an analogous result for arbitrary closed ideals $\cJ$ of $\cF$, as the following example shows.
\begin{example} \label{ex43.6}
Let $\cA := \cS(\eT(C))$ the algebra of the finite sections method for Toeplitz operators and $\cK$ the ideal of the compact sequences in the corresponding algebra $\cF$. Then
\[
\cJ := \{ (K_n) \in \cK : \lim_{n \to \infty} \|K_{2n}\| = 0 \}
\]
is a closed ideal of $\cF$. Employing again the explicit description of $\cS(\eT(C))$ in Theorem \ref{tneu} $(a)$, it is not hard to see that $\cS(\eT(C)) \cap \cJ = \cG$. Consequently, the canonical homomorphism $\pi^{\cS(\eT(C)) \cap \cJ}$ coincides with $\pi^\cG$ and is, thus, fractal. But $\cG_\eta = (\cS(\eT(C)) \cap \cJ)_\eta$ is a proper subset of $\cS(\eT(C))_\eta \cap \cJ_\eta$ for the sequence $\eta(n) := 2n+1$. Indeed, the sequence $(P_{2n+1} K P_{2n+1})$ belongs to $\cS(\eT(C))_\eta \cap \cJ_\eta$ for each compact operator $K$. \hfill \qed
\end{example}
The previous considerations suggest the following definitions. Note that both definitions coincide if $\cJ$ is a closed ideal of $\cA$ {\em and} $\cF$.
\begin{defi} \label{d43.8}
Let $\cA$ be a $C^*$-subalgebra of $\cF$. \\[1mm]
$(a)$ If $\cJ$ is a closed ideal of $\cA$ then $\cA$ is called {\em $\cJ$-fractal} \index{algebra!$\cJ$-fractal}
if the canonical homomorphism $\pi^\cJ : \cA \to \cA/\cJ$ is fractal. \\[1mm]
$(b)$ If $\cJ$ is a closed ideal of $\cF$ then $\cA$ is called {\em $\cJ$-fractal} if $\cA$ is $(\cA \cap \cJ)$-fractal and if $(\cA \cap \cJ)_\eta = \cA_\eta \cap \cJ_\eta$ for every strictly increasing sequence $\eta : \sN \to \sN$.
\end{defi}
The following results show that $\cJ$-fractality implies what one expects: A sequence in a $\cJ$-fractal algebra belongs to $\cJ$ or is invertible modulo $\cJ$ if and only if at least one of its subsequences has this property.
\begin{theo} \label{t43.9}
Let $\cJ$ be a closed ideal of $\cF$. A $C^*$-subalgebra $\cA$ of $\cF$ is $\cJ$-fractal if and only if the following implication holds for every sequence $\bA \in \cA$ and every strictly increasing sequence $\eta$
\begin{equation} \label{e43.10}
R_\eta (\bA) \in \cJ_\eta \quad \Longrightarrow \quad \bA \in \cJ.
\end{equation}
\end{theo}
{\bf Proof.} Let $\cA$ be $\cJ$-fractal and $\bA \in \cA$ a sequence with $R_\eta (\bA) \in \cJ_\eta$. Then
$R_\eta (\bA) \in \cA_\eta \cap \cJ_\eta = (\cA \cap \cJ)_\eta$, and the $(\cA \cap \cJ)$-fractality of $\cA$ implies $\bA \in \cJ$ via Theorem \ref{t43.1}.

Conversely, let (\ref{e43.10}) hold for every strictly increasing sequence $\eta$. From Theorem \ref{t43.1} we conclude that $\cA$ is $(\cA \cap \cJ)$-fractal. Further, the inclusion $\subseteq$ in $(\cA \cap \cJ)_\eta = \cA_\eta \cap \cJ_\eta$ is obvious. For the reverse inclusion, let $\bA$ be a sequence in $\cF$ with $R_\eta (\bA) \in \cA_\eta \cap \cJ_\eta$. Then there are sequences $\bB \in \cA$ and $\bJ \in \cJ$ such that $R_\eta (\bA) = R_\eta (\bB) = R_\eta (\bJ)$. Since $R_\eta (\bB) \in \cJ_\eta$, the implication (\ref{e43.10}) gives $\bB \in \cJ$. Hence, $R_\eta (\bB) \in (\cA \cap \cJ)_\eta$, and since $R_\eta (\bB) = R_\eta (\bA)$, one also has $R_\eta (\bA) \in (\cA \cap \cJ)_\eta$. \hfill \qed
\begin{theo} \label{t43.11}
Let $\cJ$ be a closed ideal of $\cF$ and $\cA$ a $\cJ$-fractal and unital $C^*$-subalgebra of $\cF$. Then the following implication holds for every sequence $\bA \in \cA$ and every strictly increasing sequence $\eta$
\begin{equation} \label{e43.12}
R_\eta (\bA) + \cJ_\eta \; \mbox{is invertible in} \; \cF_\eta/\cJ_\eta \; \Longrightarrow \; \bA + \cJ \;
\mbox{is invertible in} \; \cF/\cJ.
\end{equation}
\end{theo}
{\bf Proof.} Let $\bA \in \cA$ be such that $R_\eta (\bA)  + \cJ_\eta$ is invertible in $\cF_\eta/\cJ_\eta$. Since $C^*$-algebras are inverse closed, this coset is also invertible in $(\cA_\eta + \cJ_\eta)/\cJ_\eta$. The latter algebra is canonically $^*$-isomorphic to $\cA_\eta/(\cA_\eta \cap \cJ_\eta)$, hence, to $\cA_\eta/(\cA \cap \cJ)_\eta$ by $\cJ$-fractality of $\cA$. Thus, the coset $R_\eta (\bA) + (\cA \cap \cJ)_\eta$ is invertible in $\cA_\eta/ (\cA \cap \cJ)_\eta$. Choose sequences $\bB \in \cA$ and $\bJ \in \cA \cap \cJ$ such that
\[
R_\eta (\bA) \, R_\eta(\bB) = R_\eta(\bI) + R_\eta(\bJ)
\]
where $\bI$ denotes the identity element of $\cF$. Applying the homomorphism $\pi^{\cA \cap \cJ}_\eta$ to both sides of this equality one gets
\[
\pi^{\cA \cap \cJ} (\bA) \, \pi^{\cA \cap \cJ} (\bB) = \pi^{\cA \cap \cJ} (\bI) + \pi^{\cA \cap \cJ} (\bJ)
\]
which shows that $\bA \bB - \bI \in \cJ$. Hence, $\bA$ is invertible modulo $\cJ$ from the right-hand side. The
invertibility from the left-hand side follows analogously. \hfill \qed
\begin{coro} \label{c43.13}
Let $\cJ$ be a closed ideal of $\cF$ and $\cA$ a $\cJ$-fractal and unital $C^*$-subalgebra of $\cF$. Then a sequence $\bA \in \cA$ \\[1mm]
$(a)$ belongs to $\cJ$ if and only if there is a strictly increasing sequence $\eta$ such that $\bA_\eta$ belongs to $\cJ_\eta$. \\[1mm]
$(b)$ is invertible modulo $\cJ$ if and only there is a strictly increasing sequence $\eta$ such that $\bA_\eta$ is invertible modulo $\cJ_\eta$.
\end{coro}
We still mention the following simple facts for later reference.
\begin{prop} \label{p43.15}
Let $\cJ$ be a closed ideal of $\cF$ and $\cA$ a $\cJ$-fractal $C^*$-subalgebra of $\cF$. Then \\[1mm]
$(a)$ every $C^*$-subalgebra of $\cA$ is $\cJ$-fractal. \\[1mm]
$(b)$ if $\cI$ is an ideal of $\cF$ with $\cJ \subseteq \cI$ and if $(\cA \cap \cI)_\eta = \cA_\eta \cap \cI_\eta$ for each strictly increasing sequence $\eta : \sN \to \sN$, then $\cA$ is $\cI$-fractal.
\end{prop}
{\bf Proof.} $(a)$ Let $\cB$ be a $C^*$-subalgebra of $\cA$, and let $\bB$ be a sequence in $\cB$ with $R_\eta(\bB) \in \cJ_\eta$ for a certain strictly increasing sequence $\eta$. Then $R_\eta(\bB) \in \cA_\eta \cap \cJ_\eta$. Since $\cA$ is $\cJ$-fractal, Theorem \ref{t43.9} implies that $\bB \in \cJ$. Hence $\cB$ is $\cJ$-fractal, again by Theorem \ref{t43.9}. \\[1mm]
$(b)$ Let $R_\eta (\bA) \in \cI_\eta$ for a sequence $\bA \in \cA$ and a strictly increasing sequence $\eta$. By
hypothesis, $R_\eta (\bA) \in (\cA \cap \cI)_\eta$. Choose a sequence $\bJ \in \cA \cap \cI$ with $R_\eta (\bA) = R_\eta (\bJ)$. The $\cJ$-fractality of $\cA$ implies that $\bA - \bJ \in \cJ$, whence $\bA \in \bJ + \cJ \subseteq \cI$. By Theorem \ref{t43.9}, $\cA$ is $\cI$-fractal. \hfill \qed
\subsection{Essential fractality and Fredholm property} \label{ss43.3}
Let again $\cF$ be the algebra of matrix sequences with dimension function $\delta$ and $\cK$ the associated ideal of compact sequences. We call the $\cK$-fractal $C^*$-subalgebras of $\cF$ {\em essentially fractal}.

Note that each restriction $\cF_\eta$ of $\cF$ is again an algebra of matrix sequences (with dimension function $\delta \circ \eta$); hence, the restriction $\cK_\eta$ of $\cK$ is just the ideal of the compact sequences related with $\cF_\eta$. If we speak on {\em compact subsequences} and {\em Fredholm subsequences} in what follows, we thus mean sequences $R_\eta \bA \in \cK_\eta$ and sequences $R_\eta \bA$ which are invertible modulo $\cK_\eta$, respectively. In these terms, Corollary \ref{c43.13} reads as follows.
\begin{coro} \label{c43.16}
Let $\cA$ be an essentially fractal and unital $C^*$-subalgebra of $\cF$. Then a sequence $\bA \in \cA$ is compact (resp. Fredholm) and only if one of the subsequences of $\bA$ is compact (resp. Fredholm).
\end{coro}
The following is a consequence of Proposition \ref{p43.15}.
\begin{coro} \label{c43.16a}
Let $\cA$ be a fractal $C^*$-subalgebra of $\cF$. If $(\cA \cap \cK)_\eta = \cA_\eta \cap \cK_\eta$ for each strictly increasing sequence $\eta : \sN \to \sN$, then $\cA$ is essentially fractal.
\end{coro}
Essential fractality has striking consequences for the behavior of the smallest singular values.
\begin{theo} \label{t43.17}
Let $\cA$ be an essentially fractal and unital $C^*$-subalgebra of $\cF$. A sequence $(A_n) \in \cA$ is Fredholm if and only if there is a $k \in \sN$ such that
\begin{equation} \label{e43.18}
\limsup_{n \to \infty} \sigma_k (A_n) > 0.
\end{equation}
\end{theo}
{\bf Proof.} If $(A_n)$ is Fredholm then, by Theorem \ref{t42.1} $(c)$, $\liminf_{n \to \infty} \sigma_k \, (A_n) > 0$ for some $k \in \sN$, whence (\ref{e43.18}). Conversely, let (\ref{e43.18}) hold for some $k$. We choose a strictly increasing sequence $\eta$ such that $\lim_{n \to \infty} \, \sigma_k (A_{\eta(n)}) > 0$. Thus, the restricted sequence $(A_{\eta(n)})_{n \ge 1}$ is Fredholm by Theorem \ref{t42.1}. Since $\cA$ is essentially fractal, Corollary \ref{c43.16} $(b)$ implies the Fredholm property of the sequence $(A_n)$ itself. \hfill \qed \\[3mm]
Consequently, if a sequence $(A_n)$ in an essentially fractal and unital $C^*$-subalgebra of $\cF$ is {\em not} Fredholm, then
\begin{equation} \label{e43.19}
\lim_{n \to \infty} \, \sigma_k(A_n) = 0 \qquad \mbox{for each} \; k \in \sN.
\end{equation}
In analogy with operator theory, we call a sequence $(A_n)$ with property (\ref{e43.19}) {\em not normally solvable}. \index{sequence!not normally solvable}
\begin{coro} \label{c43.20}
Let $\cA$ be an essentially fractal and unital $C^*$-subalgebra of $\cF$. Then a sequence in $\cA$ is either Fredholm or not normally solvable.
\end{coro}
\begin{example} \label{ex43.21}
Consider the finite sections algebra $\cS(\eT(C))$ for Toeplitz operators. It is a simple consequence of Theorem \ref{tneu} $(a)$ that $(\cS(\eT(C)) \cap \cK)_\eta = \cS(\eT(C))_\eta \cap \cK_\eta$ for each strictly increasing sequence $\eta$. Since $\cS(\eT(C))$ is fractal and $\cG \subset \cK$, the algebra $\cS(\eT(C))$ is essentially fractal by Corollary \ref{c43.16a}. \hfill \qed
\end{example}
\subsection{Essential fractal restriction} \label{ss43.4}
Our final goal is an analogue of Theorem \ref{tf3.20} for essential fractality. Recall that we based the proof of Theorem \ref{tf3.20} on the fact that there is a sequence $\eta$ such that the norms $\|A_{\eta(n)}\|$ converge for each sequence $(A_n)$. We start with showing that $\eta$ can be even chosen such that not only the sequences
$(\|A_{\eta(n)}\|) = (\Sigma_1(A_{\eta(n)}))$ converge, but {\em every} sequence $(\Sigma_k(A_{\eta(n)}))$ with $k \in \sN$. Here, $\Sigma_1(A) \ge \ldots \ge \Sigma_n(A)$ denote the decreasingly ordered singular values of the $n \times n$-matrix $A$.
\begin{prop} \label{p43.22}
Let $\cA$ be a separable $C^*$-subalgebra of $\cF$. Then there is a strictly increasing sequence $\eta : \sN \to \sN$ such that the sequence $(\Sigma_k(A_{\eta(n)}))_{n \ge 1}$ converges for every sequence $(A_n)_{n \ge 1} \in \cA$ and every $k \in \sN$.
\end{prop}
{\bf Proof.} First consider a single sequence $(A_n) \in \cA$. We choose a strictly increasing sequence $\eta_1 : \sN \to \sN$ such that the sequence $(\Sigma_1(A_{\eta_1(n)}))_{n \ge 1}$ converges, then a subsequence $\eta_2$ of $\eta_1$ such that the sequence $(\Sigma_2(A_{\eta_2(n)}))_{n \ge 1}$ converges, and so on. The sequence $\eta(n) := \eta_n(n)$ has the property that the sequence $(\Sigma_k(A_{\eta(n)}))_{n \ge 1}$ converges for every $k \in \sN$.

Now let $(\bA^m)_{m \ge 1}$ be a countable dense subset of $\cA$, consisting of sequences $\bA^m = (A_n^m)_{n \ge 1}$.
We use the result of the previous step to find a strictly increasing sequence $\eta_1 : \sN \to \sN$ such that the sequences $(\Sigma_k(A^1_{\eta_1(n)}))_{n \ge 1}$ converge for every $k \in \sN$, then a subsequence $\eta_2$ of $\eta_1$ such that the sequences $(\Sigma_k(A^2_{\eta_2(n)}))_{n \ge 1}$ converge for every $k$, and so on. Then the sequence $\eta(n) := \eta_n(n)$ has the property that the sequences $(\Sigma_k(A^m_{\eta(n)}))_{n \ge 1}$ converge for every pair $k, \, m \in \sN$.

Let $\eta$ be as in the previous step, i.e., the sequences $(\Sigma_k(A^m_{\eta(n)}))_{n \ge 1}$ converge for every $k \in \sN$ and for every sequence $\bA^m = (A_n^m)_{n \ge 1}$ in a countable dense subset of $\cA$. We show that then the sequences $(\Sigma_k(A_{\eta(n)}))_{n \ge 1}$ converge for every $k \in \sN$ and every sequence $\bA = (A_n)$ in $\cA$. Fix $k \in \sN$ and let $\varepsilon > 0$. Using the well known inequality $|\Sigma_k(A) - \Sigma_k(B)| \le \|A-B\|$ we obtain
\begin{eqnarray*}
\lefteqn{|\Sigma_k(A_{\eta(n)}) - \Sigma_k(A_{\eta(l)})|} \\
&& \le |\Sigma_k(A_{\eta(n)}) - \Sigma_k(A^m_{\eta(n)})| + |\Sigma_k(A^m_{\eta(n)}) - \Sigma_k(A^m_{\eta(l)})| \\
&& \quad + \; |\Sigma_k(A^m_{\eta(l)}) - \Sigma_k(A_{\eta(l)})| \\
&& \le \|A_{\eta(n)} - A^m_{\eta(n)}\| + |\Sigma_k(A^m_{\eta(n)}) - \Sigma_k(A^m_{\eta(l)})| + \|A^m_{\eta(l)} - A_{\eta(l)}\| \\
&& \le 2 \, \|\bA - \bA^m\|_\cF + |\Sigma_k(A^m_{\eta(n)}) - \Sigma_k(A^m_{\eta(l)})|.
\end{eqnarray*}
Now choose $m \in \sN$ such that $\|\bA - \bA^m\|_\cF < \varepsilon/3$ and then $N \in \sN$ such that $|\Sigma_k(A^m_{\eta(n)}) - \Sigma_k(A^m_{\eta(l)})| < \varepsilon/3$ for all $n, \, l \ge N$. Then $|\Sigma_k(A_{\eta(n)}) - \Sigma_k(A_{\eta(l)})| < \varepsilon$ for all $n, \, l \ge N$. Thus, $(\Sigma_k(A_{\eta(n)}))_{n \ge 1}$ is a Cauchy sequence, hence convergent. \hfill \qed
\begin{prop} \label{p43.23}
Let $\cA$ be a $C^*$-subalgebra of $\cF$ with the property that the sequences $(\Sigma_k(A_n))_{n \ge 1}$ converge for every sequence $(A_n) \in \cA$ and every $k \in \sN$. Then $\cA$ is essentially fractal.
\end{prop}
{\bf Proof}. Let $\bK = (K_n) \in \cA$ and let $\eta: \sN \to \sN$ be a strictly increasing sequence such that $\bK_\eta \in \cK_\eta$. Then, by Theorem \ref{t41.12} $(b)$,
\begin{equation} \label{e43.24}
\lim_{k \to \infty} \limsup_{n \to \infty} \Sigma_k (K_{\eta(n)}) = 0.
\end{equation}
By hypothesis, $\limsup_{n \to \infty} \Sigma_k (K_{\eta(n)}) = \lim _{n \to \infty} \Sigma_k (K_n)$. Hence, (\ref{e43.24}) implies $\lim_{k \to \infty} \lim_{n \to \infty} \Sigma_k (K_n) = 0$, whence $\bK \in \cK$ by assertion $(a)$ of Theorem \ref{t41.12}. Thus, every sequence in $\cA$ which has a compact subsequence is compact itself. Thus $\cA$ is essentially fractal by Theorem \ref{t43.9}. \hfill \qed
\begin{theo} \label{t43.25}
Let $\cA$ be a separable $C^*$-subalgebra of $\cF$. Then there is a strictly increasing sequence $\eta : \sN \to \sN$
such that the restricted algebra $\cA_\eta = R_\eta \cA$ is essentially fractal.
\end{theo}
Indeed, if $\eta$ is as in Proposition \ref{p43.22}, then the restriction $\cA_\eta$ is essentially fractal by Proposition \ref{p43.23}. \hfill \qed \\[3mm]
We know from Theorems \ref{tf3.20} and \ref{t43.25} that every separable $C^*$-subalgebra of $\cF$ has both a fractal and an essentially fractal restriction. If is an open question whether this fact holds for arbitrary closed ideals $\cJ$ of $\cF$ in place of $\cG$ or $\cK$, i.e., whether one can always force $\cJ$-fractality by a suitable restriction.
\section{Essential spectral approximation} \label{s44}
In a series of papers \cite{Arv7,Arv3,Arv4}, Arveson studied the question of whether one can discover the essential spectrum of a self-adjoint operator $A$ from the behavior of the eigenvalues of the finite sections $P_n A P_n$ of $A$. More generally, one might ask whether one can discover the {\em essential spectrum} of a self-adjoint sequence $\bA = (A_n) \in \cF$ (i.e., the spectrum of the coset $\bA + \cK$, considered as an element of the quotient algebra $\cF/\cK$) from the behavior of the eigenvalues of the matrices $A_n$? To answer this question, Arveson introduced the notions of essential and transient points, and he discovered (under an additional condition) a certain dichotomy: if $A$ is a self-adjoint band-dominated operator, then every point in $\sR$ is either transient or essential; see Subsection \ref{ss41neu}. The goal of this section is to relate the essential spectral approximation with the property of essential fractality. In particular, we will see that a subalgebra $\cA$ of $\cF$ is essentially fractal if and only if every self-adjoint sequence in $\cA$ has Arveson dichotomy.
\subsection{Essential spectra of self-adjoint sequences} \label{ss44.1}
Given a self-adjoint matrix $A$ and a subset $M$ of $\sR$, let $N(A, \, M)$ \index{$N(A, \, M)$} denote the number of eigenvalues of $A$ which lie in $M$, counted with respect to their multiplicity. If $M = \{\lambda\}$ is a singleton, we write $N(A, \, \lambda)$ in place of $N(A, \, \{\lambda\})$. Thus, if $\lambda$ is an eigenvalue of $A$, then $N(A, \, \lambda)$ is its multiplicity.

Let $\bA = (A_n) \in \cF$ be a self-adjoint sequence. Following Arveson \cite{Arv7,Arv3,Arv4}, a point $\lambda \in \sR$ is called {\em essential} \index{point!essential} for this sequence if, for every open interval $U$ containing $\lambda$,
\[
\lim_{n \to \infty} N(A_n, \, U) = \infty,
\]
and $\lambda \in \sR$ is called {\em transient} \index{point!transient} for $\bA$ if there is an open interval $U$ which contains $\lambda$ such that
\[
\sup_{n \in \sN} \, N(A_n , \, U) < \infty.
\]
Thus, $\lambda \in \sR$ is {\em not} essential for $\bA$ if and only if $\lambda$ is transient for a subsequence of $\bA$, and $\lambda$ is {\em not} transient for $\bA$ if and only if $\lambda$ is essential for a subsequence of $\bA$. Moreover, if a point $\lambda$ is transient (resp. essential) for $\bA$, then is is also transient (resp. essential) for every subsequence of $\bA$.
\begin{theo} \label{t44.1}
Let $\bA \in \cF$ be a self-adjoint sequence. A point $\lambda \in \sR$ belongs to the essential spectrum of $\bA$ if and only if it is not transient for the sequence $\bA$.
\end{theo}
{\bf Proof.} Let $\bA = (A_n)$ be a bounded sequence of self-adjoint matrices. First let $\lambda \in \sR \setminus \sigma (\bA + \cK)$. We set $B_n := A_n - \lambda I_n$ and have to show that 0 is transient for the sequence $(B_n)$. Since $\lambda \in \sR \setminus \sigma (\bA + \cK)$, the sequence $(B_n)$ is Fredholm. By Theorem \ref{t42.1} $(c)$,
there is a $k \in \sZ^+$ such that
\[
\liminf_{n \to \infty} \sigma_{k+1} (B_n) =: C > 0 \quad \mbox{and} \quad \liminf_{n \to \infty} \sigma_k (B_n) = 0.
\]
Let $U := (-C/2, \, C/2)$. Since the singular values of a self-adjoint matrix are just the absolute values of the
eigenvalues of that matrix, we conclude that $N(B_n, \, U) \le k$ for all sufficiently large $n$. Thus, 0 is transient.

Conversely, let $\lambda \in \sR$ be transient for $(A_n)$. We claim that $(A_n - \lambda I_n)$ is a Fredholm
sequence. By transiency, there is an interval $U = (\lambda - \varepsilon, \, \lambda + \varepsilon)$ with $\varepsilon > 0$ such that $\sup_{n \in \sN} N(A_n, \, U) =: k < \infty$. Let $T_n$ denote the orthogonal projection from $\sC^{\delta(n)}$ onto the $U$-spectral subspace of $A_n$. Then $\rank T_n$ is not greater than $k$. It is moreover obvious that the matrices $B_n := (A_n - \lambda P_n)(I - T_n) + T_n$ are invertible for all $n \in \sN$ and that their inverses are uniformly bounded by the maximum of $1/\varepsilon$ and $1$. Hence, $(B_n^{-1}) \in \cF$ and
\begin{equation} \label{e44.1a}
(A_n - \lambda P_n)(I - T_n) B_n^{-1} = I - T_n B_n^{-1}.
\end{equation}
Since $(T_n)$ is a compact sequence (of essential rank not greater than $k$), this identity shows that the coset $(A_n - \lambda I_n) + \cK$ is invertible from the right-hand side. Since this coset is self-adjoint, it is then invertible from both sides. Thus, $(A_n - \lambda I_n)$ is a Fredholm sequence. \hfill \qed
\begin{prop} \label{p44.1a1}
The set of the non-transient points and the set of the essential points of a self-adjoint sequence $\bA \in \cF$ are compact.
\end{prop}
{\bf Proof.} The first assertion is an immediate consequence of Theorem \ref{t44.1}. The second assertion will follow once we have shown that the set of the essential points of $\bA$ is closed.

Let $(\lambda_k)$ be a sequence of essential points for $\bA = (A_n)$ with limit $\lambda$. Assume that $\lambda$ is not essential for $\bA$. Then there is a strictly increasing sequence $\eta : \sN \to \sN$ such that $\lambda$ is transient for $\bA_\eta$. Let $U$ be an open neighborhood of $\lambda$ with $\sup_{n \in \sN} \, N(A_{\eta(n)}, \, U) =: c < \infty$. Since $\lambda_k \to \lambda$ and $U$ is open, there are a $k \in \sN$ and an open neighborhood $U_k$ of $\lambda_k$ with $U_k \subseteq U$. Clearly, $N(A_{\eta(n)}, \, U_k) \le N(A_{\eta(n)}, \, U) \le c$. On the other hand, since $\lambda_k$ is also essential for the restricted sequence $\bA_\eta$, one has $N(A_{\eta(n)}, \, U_k) \to \infty$ as $n \to \infty$, a contradiction. \hfill \qed \\[3mm]
Note that the set of the non-transient points of a self-adjoint sequence is non-empty by Theorem \ref{t44.1}, whereas it is easy to construct self-adjoint sequences without any essential point: take a sequence which alternates between the zero and the identity matrix. In contrast to this observation, the following result shows that sequences which arise by discretization of a self-adjoint operator, always possess essential points. Let $H$ be an infinite dimensional separable Hilbert space with filtration $\cP := (P_n)$, and define the algebra $\cF^\cP$ as in Section \ref{s1}.
One can think of $\cF^\cP$ as a $C^*$-subalgebra of the algebra $\cF^\delta$ with dimension function $\delta (n) := \rank P_n$.
\begin{theo} \label{t44.2}
Let $\bA := (A_n) \in \cF^\cP$ be a self-adjoint sequence with strong limit $A$. Then every point in the essential spectrum of $A$ is an essential point for $\bA$.
\end{theo}
{\bf Proof.} We show that $A - \lambda I$ is a Fredholm operator if $\lambda \in \sR$ is not essential for $\bA$. Then $\lambda$ is transient for a subsequence of $\bA$, i.e., there are an infinite subset $\sM$ of $\sN$ and an interval $U = (\lambda - \varepsilon, \, \lambda + \varepsilon)$ with $\varepsilon > 0$ such that
\begin{equation} \label{e44.2a}
\sup_{n \in \sM} N(A_n, \, U) =: k < \infty.
\end{equation}
Let $T_n$ denote the orthogonal projection from $H$ onto the $U$-spectral subspace of $A_n P_n$. By (\ref{e44.2a}), the rank of the projection $T_n$ is not greater than $k$ if $n \in \sM$. So we conclude that the operators $B_n := (A_n - \lambda P_n)(I - T_n) + T_n$ are invertible for all $n \in \sM$ and that their inverses are uniformly bounded by the maximum of $1/\varepsilon$ and $1$. Hence,
\begin{equation} \label{e44.3}
(A_n - \lambda P_n)(I - T_n) B_n^{-1} = I - T_n B_n^{-1}
\end{equation}
for all $n \in \sM$. By the weak sequential compactness of the unit ball of $L(H)$, one finds weakly convergent subsequences $((I - T_{n_r}) B_{n_r}^{-1})_{r \ge 1}$ of $((I - T_n) B_n^{-1}) _{n \in \sM}$ and $(T_{n_r} B_{n_r}^{-1})_{r \ge 1}$ of $(T_n B_n^{-1})_{n \in \sM}$ with limits $B$ and $T$, respectively. The product of a weakly convergent sequence with limit $C$ and a $^*$-strongly convergent sequence with limit $D$ is weakly convergent with limit $CD$. Thus, passing to subsequences and taking the weak limit in (\ref{e44.3}) yields $(A - \lambda I) B = I - T$. Further, the rank of $T$ is not greater than $k$ by Lemma 5.7 in \cite{BGr3}. Thus, $(A - \lambda I) B - I$ is a compact operator. The compactness of $B (A - \lambda I) - I$ follows similarly. Hence, $A$ is a Fredholm operator. \hfill \qed \\[3mm]
Arveson gave a first example where the inclusion in Theorem \ref{t44.2} is proper. Specifically, he constructed a self-adjoint unitary operator $A \in L(l^2(\sN))$ with
\begin{equation}\label{e44.4}
\sigma(A) = \sigma_{ess} (A) = \{-1, \, 1\}
\end{equation}
such that 0 is an essential point of the sequence $(P_nAP_n)$.
\subsection{Arveson dichotomy and essential fractality} \label{ss41neu}
We say that a self-adjoint sequence $\bA \in \cF$ enjoys {\em Arveson's dichotomy} \index{Arveson's dichotomy} if every real number is either essential or transient for this sequence. Note that Arveson dichotomy is preserved when passing to subsequences. Arveson introduced and studied this property in \cite{Arv7,Arv3,Arv4}. In particular, he proved the dichotomy of the finite sections sequence $(P_n A P_n)$ when $A$ is a self-adjoint band-dominated operator which satisfies a Wiener and a Besov space condition. A generalization to arbitrary band-dominated operators was obtained in \cite{Roc10}.
\begin{theo} \label{t44.1a2}
The set of all self-adjoint sequences in $\cF$ with Arveson dichotomy is closed in $\cF$.
\end{theo}
{\bf Proof.} Let $(\bA_n)_{n \in \sN}$ be a sequence of self-adjoint sequences in $\cF$ with Arveson dichotomy which converges to a (necessarily self-adjoint) sequence $\bA$ in the norm of $\cF$. Then $\bA_n + \cK \to \bA + \cK$ in the norm of $\cF/\cK$. Since $\bA_n + \cK$ and $\bA + \cK$ are self-adjoint elements of $\cF/\cK$, this implies that the spectra of $\bA_n + \cK$ converge to the spectrum of $\bA + \cK$ in the Hausdorff metric. Thus, by Theorem \ref{t44.1}, the sets of the non-transient points of $\bA_n$ converge to the set of the non-transient points of $\bA$. Since the $\bA_n$ have Arveson dichotomy by hypothesis, this finally implies that the sets of the essential points of $\bA_n$ converge to the set of the non-transient points of $\bA$ in the Hausdorff metric.

Let now $\lambda$ be a non-transient point for $\bA$ and assume that $\lambda$ is not essential for $\bA$. Then there is a strictly increasing sequence $\eta : \sN \to \sN$ such that $\lambda$ is transient for the restricted sequence $\bA_\eta$. As we have seen above, there is a sequence $(\lambda_n)$, where $\lambda_n$ is an essential point for $\bA_n$, with $\lambda_n \to \lambda$. Since the property of being an essential is preserved under passage to a subsequence, $\lambda_n$ is also essential for the restricted sequence $(\bA_n)_\eta$.

Since the sequences $(\bA_n)_\eta$ also have Arveson dichotomy and since $(\bA_n)_\eta \to \bA_\eta$ in the norm of $\cF_\eta$, we can repeat the above arguments to conclude that the sets $M_n$ of the essential points for $(\bA_n)_\eta$ converge to the set $M$ of the non-transient points for $\bA_\eta$ in the Hausdorff metric. Since $\lambda_n \in M_n$ by construction, this implies that $\lambda \in M$. This means that $\lambda$ in not transient for $\bA_\eta$, a contradiction. \hfill \qed \\[3mm]
Here is the announced result which relates Arveson dichotomy with essential fractality.
\begin{theo} \label{t44.1b}
Let $\cA$ be a unital $C^*$-subalgebra of $\cF$. Then $\cA$ is essentially fractal if and only if every self-adjoint sequence in $\cA$ has Arveson dichotomy.
\end{theo}
{\bf Proof.} First let $\cA$ be essentially fractal. Let $\bA$ be a self-adjoint sequence in $\cA$ and $\lambda \in \sR$ a point which is not essential for $\bA$. The $\lambda$ is transient for a subsequence of $\bA$, thus, $0$ is transient for a subsequence of $\bA - \lambda \bI$. From Theorem \ref{t44.1} we conclude that this subsequence has the Fredholm property. Then, by Corollary \ref{c43.16} $(b)$ and since $\cA$ is essentially fractal, the sequence $\bA - \lambda \bI$ itself is a Fredholm sequence. Thus, $0$ is transient for $\bA - \lambda \bI$ by Theorem \ref{t44.1} again, whence finally follows that $\lambda$ is transient for $\bA$. Hence, $\bA$ has Arveson dichotomy.

Now assume that $\cA$ is not essentially fractal. Then, by Theorem \ref{t43.9}, there are a sequence $\bA = (A_n) \in \cA$ and a strictly increasing sequence $\eta : \sN \to \sN$ such that the restricted sequence $\bA_\eta$ belongs to $\cK_\eta$ but $\bA \not\in \cK$. The self-adjoint sequence $\bA^* \bA$ has the same properties, i.e., $(\bA^* \bA)_\eta = \bA^*_\eta \bA_\eta \in \cK_\eta$, but $\bA^* \bA \not\in \cK$.

Since $\bA^*_\eta \bA_\eta \in \cK_\eta$, the essential spectrum of $\bA^*_\eta \bA_\eta$ (i.e., the spectrum of the coset $\bA^*_\eta \bA_\eta + \cK_\eta$ in $\cF_\eta/\cK_\eta$) consists of the point $0$ only. Thus, by Theorem \ref{t44.1}, $0$ is the only non-transient point for the restricted sequence $\bA^*_\eta \bA_\eta$.

Since $\bA^* \bA \not\in \cK$, there is a strictly increasing sequence $\mu : \sN \to \sN$ such that $\mu(\sN) \cap \eta(\sN) = \emptyset$ and $\bA^*_\mu \bA_\mu \not\in \cK_\eta$. Hence, the essential spectrum of $\bA^*_\mu \bA_\mu$ contains at least one point $\lambda \neq 0$, and this point is non-transient for $\bA^*_\mu \bA_\mu$ by Theorem \ref{t44.1} again. Hence, there is a subsequence $\nu$ of $\mu$ such that $\lambda$ is essential for $\bA^*_\nu \bA_\nu$, but $\lambda \neq 0$ is transient for $\bA^*_\eta \bA_\eta$ as we have seen above. Thus, $\lambda$ is neither transient nor essential for $\bA^* \bA$. Hence, the sequence $\bA^* \bA$ does not have Arveson dichotomy. \hfill \qed
\begin{coro} \label{c44.1a3}
Every self-adjoint sequence in $\cF$ possesses a subsequence with Arveson dichotomy.
\end{coro}
{\bf Proof.} Let $\bA$ be a self-adjoint sequence in $\cF$. The smallest closed subalgebra $\cA$ of $\cF$ which contains $\bA$ is separable. By Theorem \ref{t43.25}, there is an essentially fractal restriction $\cA_\eta$ of $\cA$.  Then $\bA_\eta$ is a subsequence of $\bA$ with Arveson dichotomy by the previous theorem. \hfill \qed
{\small Author's address: \\[3mm]
Steffen Roch, Technische Universit\"at Darmstadt, Fachbereich
Mathematik, Schlossgartenstrasse 7, 64289 Darmstadt,
Germany. \\
E-mail: roch@mathematik.tu-darmstadt.de}
\end{document}